\documentclass[12pt]{article}
\usepackage{amsfonts,amsmath,amsthm,amssymb,amscd,latexsym,mathrsfs}
\usepackage[latin1]{inputenc}

\newtheorem{theorem}{Theorem}[section]
\newtheorem{proposition}{Proposition}[section]

\newtheorem{definition}{Definition}[section]
\newtheorem{example}{Example}[section]
\newtheorem{lemma}{Lemma}[section]
\begin{document}

\title{The Gauss map of Minimal graphs in the Heisenberg group}
\author{Christiam Figueroa \\\\Pontificia Universidad Católica del Perú \\
e-mail: cfiguer@pucp.edu.pe}

\maketitle

\begin{abstract}
In this paper we study some geometric properties of surfaces in the
Heisenberg group, $\mathcal{H}_{3}.$ We obtain, using the Gauss map for Lie
groups,   a partial classification of  minimal graphs in $\mathcal{H}_{3}.$
We also proof the non existence of umbilical surfaces in $\mathcal{H}_{3}.$

\end{abstract}

\section{Introduction}

The clasical Heisenberg group, $\mathcal{H}_{3}$, is the group of $3\times3$ matrices of the form
$$\begin{bmatrix}
    1 & x & z \\
    0 & 1 & y \\
    0 & 0 & 1 \\
  \end{bmatrix}
  $$
  which is diffeomorphic to $\mathbb{R}^{3}$. Besides physical reason,this group is the natural generalization of the Euclidean space, $\mathbb{R}^{3}$
   for several reasons, among which:
   It is a 2-step nilpotent (or quasi-abelian) Lie group, which is the nearest condition to be abelian. And for every left invariant
   Riemannian metric, the isometry group is 4-dimensional, maximal dimension for the isometry group of a non constant curvature metric
    in $\mathbb{R}^{3}. $  It´s know that there is no 3-dimensional Riemannian manifold with isometry group of dimension 5, so $\mathcal{H}_{3}$
     has isometry group of the largest possible dimension for a non-constant curvature space. Actually, this one of the eight Thurston geometries.

In this paper we will fix a left invariant Riemannian metric in $\mathcal{H}_{3}$ and study the geometry of surfaces with special emphasis on minimal
surfaces and the relationship with their Gauss map.

 We have organized the paper as follows. Section 2 we present the basic geometry of the Heisenberg group, $\mathcal{H}_{3}$ including
 a basis for left invariant fields.

 In section 3 we study the non parametric surfaces in $\mathcal{H}_{3}$. We calculate the coefficients of the first and second fundamental form and
 the gaussian curvature of this type of surface.

 In section 4 we present the the Gauss map for hypersurface of any Lie group and present a relationship between this map and the second fundamental form
 and give a direct proof of a non existence of umbilical surface in $\mathcal{H}_{3}$.

 In section 5, we present the classification of minimal graphs in $\mathcal{H}_{3}$ when the rank of its Gauss map is zero and one and finally we present some conditions in order that a minimal graph is a plane.
 \section{The Geometry of the Heisenberg group}
 The 3-dimensional Heisenberg group $\mathcal{H}_{3}$ a two-step nilpotent Lie group, has the following standard representation
 in $GL_{3}(\mathbb{R})$
 $$\left[
      \begin{array}{ccc}
        1 & r & t \\
        0 & 1 & s \\
        0 & 0 & 1 \\
      \end{array}
    \right]
$$
with $r,s,t\in \mathbb{R}$.

In order to describe a left-invariant metric on $\mathcal{H}_{3}$, we note that the Lie algebra $\mathfrak{h}_{3}$ of $\mathcal{H}_{3}$ is given by the
matrices
$$A=\left[
      \begin{array}{ccc}
        0 & x & z \\
        0 & 0 & y \\
        0 & 0 & 0 \\
      \end{array}
    \right] $$
    with $x,y,z$ real. The exponential map $exp:\mathfrak{h}_{3}\rightarrow \mathcal{H}_{3}$ is a global diffeomorphism, and is given by
    $$exp(A)=I+A+\frac{A^{2}}{2}=\left[
                                   \begin{array}{ccc}
                                     1 & x & z+\frac{xy}{2} \\
                                     0 & 1 & y \\
                                     0 & 0 & 1 \\
                                   \end{array}
                                 \right].$$
    Using the exponential map as a global parametrization, with the identification of the Lie algebra $\mathfrak{h}_{3}$ with $\mathbb{R}^{3}$ given by
$$(x,y,z)\leftrightarrow \left[
      \begin{array}{ccc}
        0 & x & z \\
        0 & 0 & y \\
        0 & 0 & 0 \\
      \end{array}
    \right],$$
    the group structure of $\mathcal{H}_{3}$ is given by
\begin{equation}\label{pr}
(x_{1},y_{1},z_{1})\ast (x_{2},y_{2},z_{2})=(x_{1}+x_{2},y_{1}+y_{2},z_{1}+z_{2}+\frac{x_{1}y_{2}-x_{2}y_{1}}{2}).
\end{equation}
    From now on, modulo the identification given by $exp$, we consider $\mathcal{H}_{3}$ as $\mathbb{R}^{3}$ with the product given in (\ref{pr}).
    Notice, in this model, the 1-parameter subgroups are straight lines. The Lie   algebra bracket, in terms of the canonical basis $\{e_{1},e_{2},e_{3}\}$ of $\mathbb{R}^{3}$, is given by
    $$[e_{1},e_{2}]= e_{3},  [e_{i},e_{3}]=0$$
    with $i=1,2,3.$ Now, using $\{e_{1},e_{2},e_{3}\}$ as the orthonormal frame at the identity, we have the following left-invariant metric $ds^{2}$ in $\mathcal{H}_{3},$
$$ds^{2}=dx^{2}+dy^{2}+(\frac{1}{2}ydx-\frac{1}{2}xdy+dz)^{2}.$$
And the basis of the orthonormal left-invariant vector fields is given by
    $$E_{1}=\frac{\partial}{\partial x}-\frac{y}{2}\frac{\partial}{\partial z},  E_{2}=\frac{\partial}{\partial x}+\frac{x}{2}\frac{\partial}{\partial z},
    E_{3}=\frac{\partial}{\partial z},$$
Then the Riemann connection of $ds^{2}$, in terms of the basis $\{E_{i}\}$, is given by:
$$\begin{array}{ccccc}
  \nabla_{E_{1}}E_{2} & = & \frac{1}{2}E_{3} & = & -\nabla_{E_{2}}E_{1} \\
  \nabla_{E_{1}}E_{3} & = &  -\frac{1}{2}E_{2}& = & \nabla_{E_{3}}E_{1} \\
  \nabla_{E_{2}}E_{1} & = & \frac{1}{2}E_{1} & = &\nabla_{E_{3}}E_{2}
\end{array}
$$
and $\nabla_{E_{i}}E_{i}=0$ for $i=1,2,3$.

\section{Non parametric surface in $\mathcal{H}_{3}$}
Let $S$ be a graph of a smooth function $f:\Omega \rightarrow \mathbb{R}$ where $\Omega $ is an open set of $\mathbb{R}^{2}$. We consider the
following parametrization of $S,$
\begin{equation}\label{paramet}
 X\left( x,y\right) =( x,y,f( x,y),(x,y))\in \Omega.
\end{equation}
A basis of the tangent space $T_{p}S$ associated to this
parametrization is given by
\begin{equation}
\begin{array}{ccccc}
X_{x} & = & \left( 1,0,f_{x}\right) & = & E_{1}+\left( f_{x}+\frac{y}{2}%
\right) E_{3} \\
X_{y} & = & \left( 0,1,f_{y}\right) & = & E_{2}+\left( f_{y}-\frac{x}{2}%
\right) E_{3},%
\end{array}
\label{basis}
\end{equation}%
\noindent
and
\begin{equation}
\eta \left( x,y\right) =-\left( \frac{f_{x}+\frac{y}{2}}{w}\right)
E_{1}-\left( \frac{f_{y}-\frac{x}{2}}{w}\right) E_{2}+\frac{1}{w}E_{3}
\label{normal}
\end{equation}%
where $w=\sqrt{1+\left( f_{x}+\frac{y}{2}\right) ^{2}+\left( f_{y}-\frac{x}{2%
}\right) ^{2}},$ is the unit normal vector field. Then the
coefficients of the first fundamental form of $S$  are given by%
\begin{equation}
\begin{array}{ccccl}
E & = & <X_{x},X_{x}> & = & 1+\left( f_{x}+\frac{y}{2}\right) ^{2} \\
F & = & <X_{y},X_{x}> & = & \left( f_{x}+\frac{y}{2}\right) \left( f_{y}-%
\frac{x}{2}\right) \\
G & = & <X_{y},X_{y}> & = & 1+\left( f_{y}-\frac{x}{2}\right) ^{2}.%
\end{array}
\label{1ffund}
\end{equation}%
If $\nabla $ is the Riemannian connection of $\left( \mathcal{H}%
_{3},ds^{2}\right) $, by Weingarten´s  formula for hypersurfaces, we have that%
\[
A_{\eta }v=-\nabla _{v}\eta ,\ \ \ \ v\in T_{p}S
\]%
and the coefficients of the second fundamental forma are given by
\begin{equation}
\begin{array}{ccccc}
L & = & -<\nabla _{X_{x}}\eta ,X_{x}> & = & \frac{f_{xx}+( f_{y}-\frac{x%
}{2})( f_{x}+\frac{y}{2})}{w} \\\\
M & = & -<\nabla _{X_{x}}\eta ,X_{y}> & = & \frac{f_{xy}+\frac{1}{2}\left(
f_{y}-\frac{x}{2}\right) ^{2}-\frac{1}{2}\left( f_{x}+\frac{y}{2}\right) ^{2}%
}{w} \\\\
N & = & -<\nabla _{X_{y}}\eta ,X_{y}> & = & \frac{f_{yy}-\left( f_{y}-\frac{x%
}{2}\right) \left( f_{x}+\frac{y}{2}\right) }{w}.%
\end{array}
\label{2ffund}
\end{equation}

To end this section we calculate the sectional curvature of a nonparametric surface in $\mathcal{H}_{3}$. You can see the same formula in \cite{bekkar2}.
\begin{theorem}
Let $S$ be a nonparametric surface in $\mathcal{H}_{3}$ given by $(x,y,f(x,y))$ with $(x,y)\in \Omega\subset \mathbb{R}^{2}.$ Then the sectional curvature of $S$ is given by
\begin{align*}
  w^{4}K & =w^{2}(f_{xy}^{2}-f_{xx}f_{yy}-\frac{1}{4})-(1+q^{2}){(f_{xy}+\frac{1}{2})^{2}-f_{xx}f_{yy}} \\
   & -(1+p^{2}){(f_{xy}-\frac{1}{2})^{2}-f_{xx}f_{yy}}+pq(f_{yy}-f_{xx})
\end{align*}

\noindent where $p,q$ and $w$ are defined by
$$
p=f_{x}+\frac{y}{2}, q=f_{y}-\frac{x}{2}, w=\sqrt{1+p^{2}+q^{2}}.
$$
\end{theorem}
\begin{proof}
We recall the following formula from the Gauss equation for
isometric immersions for this case
$$ K(X_{x},X_{y})-\overline{K}(X_{x},X_{y})=detA_{\eta}$$
where ${X_{x},X_{y}}$ is the basis  of  $S$, associated to the parametrization (\ref{paramet}),
$K$ and $\overline{K}$ are  the sectional curvatures of $S$ and $\mathcal{H}_{3}$ respectively.
Using this basis,  we have
$$\nabla_{X_{x}}(\nabla_{X_{x}}X_{x})=[\frac{1}{2}f_{xx}-\frac{1}{2}(f_{y}-\frac{x}{2})(f_{x}+\frac{y}{2})]E_{1}
-(f_{xy}+\frac{1}{2})E_{2}+f_{xxy}E_{3}.$$
In the same way
\begin{align*}
  \nabla_{X_{x}}(\nabla_{X_{y}}X_{x})& = [\frac{1}{2}f_{xx}-\frac{1}{4}(f_{y}-\frac{x}{2})(f_{x}+\frac{y}{2})]E_{1}
  - [f_{xy}+\frac{1}{4}(f_{x}+\frac{y}{2})^{2}-\frac{1}{4}]E_{2} \\
    & + [f_{xxy}-\frac{1}{4}(f_{y}-\frac{x}{2})]E_{3}
\end{align*}
\noindent On the other hand, $[X_{x},X_{y}]=0$. So the curvature tensor of $\mathcal{H}_{3}$, is
$$R(X_{x},X_{y})X_{x}=-\frac{1}{4}(f_{y}-\frac{x}{2})(f_{x}+\frac{y}{2})E_{1}
+[\frac{1}{4}(f_{x}+\frac{y}{2})^{2}-\frac{3}{4}]E_{2}+\frac{1}{4}(f_{y}-\frac{x}{2})E_{3},$$
so,  the sectional curvature of $\mathcal{H}_{3}$ is given by
$$\overline{K}(X_{x},X_{y})=\frac{<R(X_{x},X_{y})X_{x},X_{y}>}{\|X_{x}\wedge X_{y}\|^{2}}=\frac{1}{4}-\frac{1}{w^{2}}.$$
On the other hand, using (\ref{1ffund}) and (\ref{2ffund}), we have
$$detA_{\eta}=\frac{LN-M^{2}}{EG-F^{2}}=\frac{f_{xx}f_{yy}+pq(f_{yy}-f_{xx})-
\frac{1}{4}(p^{2}+q^{2})^{2}-f_{xy}^{2}-f_{xy}(q^{2}-p^{2})}{w^{4}}.$$
So the sectional curvature of $S$, satisfy
$$w^{4}K=f_{xx}f_{yy}-f_{xy}^{2}+\frac{1}{4}+pq(f_{yy}-f_{xx})+p^{2}(f_{xy}-\frac{1}{2})-q^{2}(f_{xy}+\frac{1}{2})-1.$$
From this relation  follows the formula.
\end{proof}
Using the above formula, J. Inoguchi classified  flat translation invariant surfaces and F. Dillen, J. Van der Veken, constructed some examples of semi-parallel surfaces in $\mathcal{H_{3}}$, see \cite{inog} and \cite{frdi} respectively.

\section{The Gauss map}

Recall the Gauss map is a function  from an oriented surface, $S\subset \mathbf{R}^{3}$ , to the unit sphere in the Euclidean space . It associates to every point on the surface its oriented unit normal vector. Considering the Euclidean space as a commutative Lie group, the Gauss map is
just the translation of the unit normal vector at any point of the surface to the origin, the identity element of $\mathbf{R}^{3}$. Reasoning in this way we define a Gauss map in the following form:

\begin{definition}
Let $S\subset G$ be an orientable hypersurface of a n-dimensional Lie group $%
G,$ provided with a left invariant metric. The map%
\[
\gamma =S\rightarrow S^{n-1}=\left\{ v\in \tilde{g}:\left\vert v\right\vert
=1\right\}
\]%
where $\gamma \left( p\right) =dL_{p}^{-1}\circ \eta \left( p\right) $, $\tilde{g}$ the Lie algebra of $G$
and $\eta $ the unitary normal vector field of $S,$ is called the Gauss map of $S.$
\end{definition}

\noindent We observe that
$$d\gamma\left( T_{p}S\right) \subseteq T_{\gamma \left( p\right) }S^{n-1}
 =  \left\{ \gamma \left( p\right) \right\} ^{\perp }  =
dL_{p}^{-1}\left( T_{p}S\right),$$

\noindent therefore $dL_{p}\circ d\gamma \left( T_{p}S\right)
\subseteq T_{p}S$ .

We know that in the Euclidean case the differential of the Gauss
map is just the second fundamental form for surfaces in
$\mathbb{R}^{3},$ this fact can be generalized for hypersurfaces
in any Lie group. The following theorem, see \cite{ripo},
states a relationship between the Gauss map and the extrinsic
geometry of $S.$

\begin{theorem}
\label{gauss}Let $S$ be an orientable hypersurfaces of a Lie group. Then
$$
dL_{p}\circ d\gamma _{p}\left( v\right) =-\left( A_{\eta }\left( v\right)
+\alpha _{\bar{\eta}}\left( v\right) \right) ,\ \ \ v\in T_{p}S,
$$
where $A_{\eta }$ is the Weingarten operator, $\alpha _{\bar{\eta}}\left(
v\right) =\nabla _{v}\bar{\eta}$ and $\bar{\eta}$ is the left invariant
vector field such that $\eta \left( p\right) =\bar{\eta}\left( p\right).$
\end{theorem}

In the case of orientable surfaces in $\mathcal{H}_{3}$ we shall obtain the expressions of
the operators $dL_{p}\circ d\gamma _{p}$ and $a_{\bar{\eta}}$ , when such a surface
is the graph of a smooth function $f\left(x,y\right) .$ In fact, using the basis $\left\{ X_{x},X_{y}\right\} ,$
given by the parametrization $\left( \ref{paramet}\right) ,$ we have%
\[
\begin{array}{ccc}
d\gamma _{p}\left( X_{x}\right) & = & \sum\limits_{i=1}^{3}\frac{\partial
a_{i}}{\partial x}E_{i}\left( e\right) \\
d\gamma _{p}\left( X_{y}\right) & = & \sum\limits_{i=1}^{3}\frac{\partial
a_{i}}{\partial y}E_{i}\left( e\right)%
\end{array}%
\]%
where $a_{i}$ are the components of the normal $\eta ,$ see $\left( \ref%
{normal}\right) ,$ and $p\in S$. Hence%
\[
\begin{array}{ccc}
dL_{p}\circ d\gamma _{p}\left( X_{x}\right) & = & \sum\limits_{i=1}^{3}\frac{%
\partial a_{i}}{\partial x}E_{i}\left( p\right) \\
dL_{p}\circ d\gamma _{p}\left( X_{y}\right) & = & \sum\limits_{i=1}^{3}\frac{%
\partial a_{i}}{\partial y}E_{i}\left( p\right).%
\end{array}%
\]%
On the other hand, we have that $dL_{p}\circ d\gamma \left( T_{p}S\right)
\subseteq T_{p}S,$ so
\[
\begin{array}{ccc}
dL_{p}\circ d\gamma _{p}\left( X_{x}\right) & = & aX_{x}+bX_{y} \\
dL_{p}\circ d\gamma _{p}\left( X_{y}\right) & = & cX_{x}+dX_{y}%
\end{array}%
\]%
Using $\left( \ref{basis}\right) $ and comparing the above two systems, we
obtain the following matrix%
\[
dL_{p}\circ d\gamma _{p}=\left[
\begin{array}{cc}
-\left( \frac{f_{x}+\frac{y}{2}}{w}\right) _{x} & -\left( \frac{f_{x}+\frac{y%
}{2}}{w}\right) _{y} \\
-\left( \frac{f_{y}-\frac{x}{2}}{w}\right) _{x} & -\left( \frac{f_{y}-\frac{x%
}{2}}{w}\right) _{y}%
\end{array}%
\right]
\]%
Notice that
\begin{equation}
\det \left( dL_{p}\circ d\gamma _{p}\right) =\frac{f_{xx}f_{yy}-f_{xy}^{2}+%
\frac{1}{4}}{w^{2}}  \label{rank}
\end{equation}%
and we will call this expression, the rank of the Gauss map. And $\alpha _{\bar{\eta}}$ is given
by the following matrix,
\[
\alpha _{\bar{\eta}}=\frac{1}{2w}\left[
\begin{array}{lr}
-\left( f_{x}+\frac{y}{2}\right) \left( f_{y}-\frac{x}{2}\right) & 1-\left(
f_{y}-\frac{x}{2}\right) ^{2} \\
\left( f_{x}+\frac{y}{2}\right) ^{2}-1 & \left( f_{x}+\frac{y}{2}\right)
\left( f_{y}-\frac{x}{2}\right)%
\end{array}%
\right]
\]%
\label{alfa}where $w$ is like ( \ref{normal}). Observe that in the
case of the Heisenberg group the trace of $\alpha _{\bar{\eta}}$
is zero. Our first result is

\begin{theorem}
\label{vertical}The vertical plane is the unique connected surface in $%
\mathcal{H}_{3}$ with the property that its Gauss map is constant.
\end{theorem}

\begin{proof}
Let $S$ be a surface in $\mathcal{H}_{3}$ parameterized as the graph of a smooth function $f\left( x,y\right) .$
As we have seen, a basis of the tangent space of $S$ is given by:%
\[
\begin{array}{ccc}
X_{x} & = & E_{1}+\left( f_{x}+\frac{y}{2}\right) E_{3} \\
X_{y} & = & E_{2}+\left( f_{y}-\frac{x}{2}\right) E_{3}.%
\end{array}%
\]%
Now if there is $p\in S$ such that $d\gamma _{p}=0,$ then $dL_{p}^{-1}\left(
T_{p}S\right) $ is a subalgebra of $h_{3},$ see \cite{ripo}. But this is a
contradiction, because $\left[ dL_{p}^{-1}\left( X_{x}\right)
,dL_{p}^{-1}\left( X_{y}\right) \right] =e_{3}\notin dL_{p}^{-1}\left(
T_{p}S\right) .$ Therefore, there are no graphs in $\mathcal{H}_{3}$ such
that its Gauss map is constant.

Now we consider $S$ as a vertical surface. In this case we can consider such
a surface as a ruled surface. We parameterize the surface by%
$$
X\left( t,s\right) =\left( t,a\left( t\right) ,s\right) ,\ \ \ \left(
t,s\right) \in U,
$$
where $U$ is an open set of $\mathbb{R}^{2}$. So the basis associated to this parametrization is
\[
\begin{array}{ccl}
X_{t} & = & E_{1}+\dot{a}E_{2}+\left( a-t\dot{a}\right) E_{3} \\
X_{s} & = & E_{3}, \\
\end{array}%
\]%
and the normal field to this surface is
$$\eta  = \frac{\dot{a}}{\sqrt{1+\left( \dot{a}\right) ^{2}}}E_{1}-\frac{1}{%
\sqrt{1+\left( \dot{a}\right) ^{2}}}E_{2}.
$$
Notice that $\eta $ is constant iff $\dot{a}\left( t\right) $ is constant,
that is, $a\left( t\right) $ is affine.
\end{proof}

To end this section, we have an alternative proof of the following result.

\begin{theorem}
There are no totally umbilical surfaces in $\mathcal{H}_{3}$
\end{theorem}

\begin{proof}
Let $S$ be an umbilical surface which is, locally, the graph of a
differentiable function $f.$ Then
\[
\begin{array}{ccc}
A_{\eta }\left( X_{x}\right) & = & \lambda X_{x} \\
A_{\eta }\left( X_{y}\right) & = & \lambda X_{y},%
\end{array}%
\]%
where $\lambda $ is a differentiable function. The Codazzi equation (in the
umbilical case) is given by%
\[
R\left( X_{x},X_{y}\right) \eta =X_{y}\left( \lambda \right)
X_{x}-X_{x}\left( \lambda \right) X_{y}=\lambda _{y}X_{x}-\lambda _{x}X_{y}.
\]%
By replacing $\left( \ref{basis}\right) $ and $\left( \ref{normal}\right) $
into the last expression, we obtain%
\[
\begin{array}{ccc}
\lambda _{x} & = & -\left( f_{x}+\frac{y}{2}\right) /w \\
\lambda _{y} & = & -\left( f_{y}-\frac{x}{2}\right) /w.%
\end{array}%
\]%
But this implies that $dL_{p}\circ d\gamma _{p}$ has a symmetric matrix in
the basis $\left( \ref{basis}\right) $. Then using theorem ( \ref%
{gauss}), we conclude that $\alpha _{\bar{\eta}}$ is also
symmetric.
Therefore, using $\left( \ref{alfa}\right) $,%
\[
w=\sqrt{3}
\]%
This yields, together with the fact that $\lambda _{xy}=\lambda _{yx},$%
\[
\left( f_{x}+\frac{y}{2}\right) _{y}=\left( f_{y}-\frac{x}{2}\right) _{x}
\]%
which is a contradiction.

Now assume that $S$ is a vertical surface, we may consider it as a ruled
surface, where the vertical lines are the rulings and the directrix, $%
a\left( t\right),$ lies in the $xy-plane.$ As usual, we parameterize the
surface by%
\begin{equation}
X\left( t,s\right) =\left( t,a\left( t\right) ,s\right) ,\ \ \ \left(
t,s\right) \in U\subseteq \mathbb{R}^{2}.  \label{verti}
\end{equation}

\noindent The coefficients of the first fundamental form in the basis $\left\{
X_{t},X_{s}\right\} $ are given by:

\begin{equation}\label{vprimfun}
\begin{array}{ccl}
E & = & 1+\dot{a}^{2}+\left( a-t\dot{a}\right) ^{2}/4 \\
F & = & \left( a-t\dot{a}\right) /2 \\
G & = & 1%
\end{array}%
\end{equation}

and the coefficients of the second fundamental forma in the same basis are
given by%
\begin{equation}\label{vsegfun}
\begin{array}{ccl}
L & = & (( a-t\dot{a})( 1+\dot{a}^{2}) -2\ddot{a})/2%
\sqrt{1+\dot{a}^{2}} \\
M & = & \sqrt{1+\dot{a}^{2}}/2 \\
N & = & 0.%
\end{array}%
\end{equation}
Since we have assumed that the surface was umbilical, we conclude that the
Weingarten operator is a diagonal matrix in any basis, in particular in the
basis associated to the parametrization $\left( \ref{verti}\right),$ then $%
0=NF-MG$ , but this is a contradiction.
\end{proof}

We remark that P. Piu proved that there are no totally geodesic hypersurfaces in $\mathcal{H}_{2n+1},$ see \cite{piu} and A. Sanini  generalized this result, that is, there are no totally umbilical hypersurfaces in this group, see \cite{sani}. Finally J. Van der Veken gave a full local classification of totally umbilical surfaces in 3-dimensional homogeneous spaces with 4-dimensional group, see \cite{jovan}. Also, R. Souam and E. Toubiana , obtained the same result independently, see \cite{suan}.

\section{Minimal graphs in $\mathcal{H}_{3}$}

We recall firstly the mean curvature formula of any surface of $\mathcal{H}_{3}$ in terms of the coefficients of
their first and second fundamental forms,
\begin{equation}
H=\frac{1}{2}\left( \frac{EN+GL-2FM}{EG-F^{2}}\right).  \label{mean curvatue}
\end{equation}
If the surface is the graph of a smooth function, using $\left( \ref{1ffund}\right) $ and $\left(\ref{2ffund}\right) $
into the above equation, we obtain the equation of the minimal graphs in $%
\mathcal{H}_{3}$
\begin{equation}
\left( 1+\left( f_{y}-\frac{x}{2}\right) ^{2}\right) f_{xx}-2\left( f_{y}-%
\frac{x}{2}\right) \left( f_{x}+\frac{y}{2}\right) f_{xy}+\left( 1+\left(
f_{x}+\frac{y}{2}\right) ^{2}\right) f_{yy}=0.  \label{minec}
\end{equation}
Before presenting some consequences of the above equation, we shall show
some examples of minimal graphs and compute the rank of their Gauss map,
using formula $\left( \ref{rank}\right). $

\begin{example}
As in Euclidean space $\mathbb{E}^{3},$ the planes $f\left( x,y\right)
=ax+by+c$ are minimal graphs of $\mathcal{H}_{3}.$ The rank of the Gauss map
is $2.$
\end{example}

Another minimal graphs may be obtained by searching for solutions of Scherk
type, i.e for solutions of the form $f\left( x,y\right) =u\left( x\right)
+v\left( y\right) +\frac{xy}{2}.$ From this method we find, among others,
the following example, see \cite{bekkar}.

\begin{example}
A surface of saddle type:
$$f\left( x,y\right) =\frac{xy}{2}+k\left[ \ln
\left( y+\sqrt{1+y^{2}}\right) +y\sqrt{1+y^{2}}\right] ,$$ where
$k\in \mathbb{R}$. Notice that this minimal surface is ruled by
affine lines, i.e. translations of 1-parameter subgroups. The rank
of its Gauss map is 1.
\end{example}

Unlike the case of Euclidean spaces, where the only complete
minimal graphs are linear (Bernstein´s  theorem), we have several
solutions defined on the entire $xy$-plane.

Let us come back to the minimal graphs equation $\left( \ref{minec}\right) .$
Notice that this is a quasilinear, elliptic P.D.E. with analytic coefficients and
therefore its solutions are analytic and satisfy the following maximum
principle, see \cite{jost}.

\begin{theorem}
Consider an elliptic, differential equation of the form%
$$
F[u]:=F\left( x,y,u,Du,D^{2}u\right)  = 0
$$
with $F:S=\Omega\times \mathbb{R}\times\mathbb{R}^{2}\times S(2,\mathbb{R})\rightarrow \mathbb{R}$
where $S(2,\mathbb{R})$ is the space of symmetric, real-valued, $2\times2$ matrices.
Let $u_{0},u_{1}\in C^{2}(\Omega)\cap C^{0}(\overline{\Omega})$, and suppose
\begin{enumerate}
  \item $F\in C^{1}(S),$
  \item $F$ is elliptic at all functions $tu_{1}+(1-t)u_{0}, 0\leq t\leq1,$
  \item $\frac{\partial F}{\partial u}\leq 0$ in $\Omega .$
\end{enumerate}
If $u_{1}\leq u_{0}$ on $\partial(\Omega)$ and $F[u_{1}]\geq F[u_{0}]$ in $\Omega$,
then either $u_{1}<u_{0}$ in $\Omega$ or $u_{0}\equiv u_{1}$ in $\Omega.$
\end{theorem}

\noindent A simple application of these facts imply the following results

\begin{theorem}
There are no compact (i.e. bounded and closed) minimal surfaces in $\mathcal{H}_{3}$
\end{theorem}

\begin{proof}
Suppose $S$ is a compact minimal surface (without boundary) in $\mathcal{H}%
_{3}.$ Take the plane $z=c,$ which is a minimal surface, such that the plane
is tangent to $S$ and $S$ lies below the plane, so by the maximum principle,
$S$ locally coincides with the plane and, by analyticity, $S$ is the plane,
which contradicts compactness.
\end{proof}

Unlike the minimal surface case, for graphs of non-zero constant mean
curvature, we have a Bernstein type theorem, see \cite{figuer}.

\begin{theorem}
There are not complete graphs of constant mean curvature $H\neq 0.$
\end{theorem}

We shall now study the stability of minimal graphs. To explain this, we need
to characterize the minimal surfaces as solution of a variational problem.
Let $S$ be a surface given by $z=f\left( x,y\right) $ with $(x,y)\in \Omega\subset \mathbb{R}^{2}.$ Then we
consider the following variation of $S:$%
\[
S_{t}\left( x,y\right) =\left( x,y,f\left( x,y\right) +th\left( x,y\right)
\right) ,\ \ \ \ \left( x,y\right) \in \Omega
\]%
where $h\in C^{1}$ and $\left. h\right\vert _{\partial \Omega }=0.$
Furthermore, the area of $S_{t}$ over $\bar{\Omega}$ is
\[
A\left( t\right) =\iint\limits_{\bar{\Omega}}w\left( t\right) dxdy
\]%
where $w\left( t\right) =\sqrt{1+\left( f_{x}+th_{x}+\frac{y}{2}\right)
^{2}+\left( f_{y}+th_{y}-\frac{x}{2}\right) ^{2}}.$ Since $S$ has least area
among all surfaces of $S_{t}$ , we have that $S$ must be critical point of $%
A\left( t\right) $ i.e. $A^{\prime }\left( 0\right) =0.$ Now we compute the
first derivative of $A\left( t\right) $%
\begin{equation}
A^{\prime }\left( t\right) =\iint\limits_{\bar{\Omega}}w^{-1}\left[ \left(
f_{x}+th_{x}+\frac{y}{2}\right) h_{x}+\left( f_{y}+th_{y}-\frac{x}{2}\right)
h_{y}\right] dxdy  \label{variacion}
\end{equation}%
Evaluated at $t=0$, integrating by parts, and using the fact that
$h=0$ on $\partial \Omega$, we find

$$
\begin{array}{ccl}
  A^{\prime }(0) & = &\iint\limits_{\bar{\Omega}}w^{-1}\left[ \left(
f_{x}+\frac{y}{2}\right) h_{x}+\left( f_{y}-\frac{x}{2}\right)
h_{y}\right] dxdy \\\\
    & =&\iint\limits_{\bar{\Omega}}\left[ (\frac{
f_{x}+\frac{y}{2}}{w})_{x}+(\frac{ f_{y}-\frac{x}{2}}{w})_{y}
\right]h dxdy.
\end{array}
$$

\noindent It follows that the equation
$$
(\frac{f_{x}+\frac{y}{2}}{w})_{x}+(\frac{
f_{y}-\frac{x}{2}}{w})_{y}=0
$$
must hold $\forall (x,y)\in \Omega $ and yields the same equation in
(\ref{minec}). Using this equation and the matrix representation of $dL_{p}\circ d\gamma_{p}$ we conclude the following

\begin{proposition}
Let $ f:\Omega\rightarrow \mathbb{R}$ be a smooth function. A graph of $f$ is a minimal surface in $ \mathcal{H}_{3}$  if and only if
the trace of $ dL_{p}\circ d\gamma_{p}$  is equal to zero.
\end{proposition}

 Now we are ready to prove the following

\begin{proposition}
Every minimal graph in $\mathcal{H}_{3}$ is stable.
\end{proposition}

\begin{proof}
It is sufficient to consider the second derivative of the area function, $%
A\left( t\right) ,$ evaluated at $t=0.$ From $\left( \ref{variacion}\right) $
we obtain that
\[
A^{\prime \prime }\left( 0\right) =\iint\limits_{\bar{\Omega}}\frac{%
h_{x}^{2}+h_{y}^{2}+\left( \left( f_{y}-\frac{x}{2}\right) h_{x}-\left(
f_{x}+\frac{y}{2}\right) h_{y}\right) ^{2}}{w^{2}}dxdy
\]%
Notice that $A^{\prime \prime }\left( 0\right) \geq 0$ and is equal to $0$ $%
iff$ $h_{x}=h_{y}=0$ that is $h=0,$ because $\left. h\right\vert _{\partial
\Omega }=0$
\end{proof}

Now we will give a classification of the minimal surfaces in $\mathcal{H}%
_{3} $ with Gauss map of rank 0 and 1. We begin with the rank 0 case

\begin{theorem}
The vertical plane is the only minimal surface in $\mathcal{H}_{3}$ with the
property that its Gauss map is constant.
\end{theorem}

\begin{proof}
As we have seen in theorem $\left( \ref{vertical}\right) $ the vertical
plane is the unique connected surface in $\mathcal{H}_{3}$ with the property
that its Gauss map is constant. Then it remains to prove that such a surface
is minimal. To do this we parameterize the surface by
$$
X\left( t,s\right) =\left( t,kt,s\right)
$$
with $k\in\mathbb{R}$ and compute, using (\ref{vprimfun}) and (\ref{vsegfun}),  the coefficients of the first and second fundamental forms in
the basis $\left\{ X_{t},X_{s}\right\} ,$%
\[
\begin{array}{cclclcccccc}
E & = & 1+k^{2}+\frac{k^{2}t^{2}}{4} & , & F & = & -\frac{kt}{2} & , & G & =
& 1 \\\\
L & = & \frac{-kt\sqrt{1+k^{2}}}{2} & , & M & = & \frac{\sqrt{1+k^{2}}}{2} &
, & N & = & 0%
\end{array}%
\]%
Finally, substituting these into the mean curvature formula
( \ref{mean curvatue}) , we have that $H=0.$
\end{proof}
Actually, the only minimal vertical surface in $\mathcal{H}_{3}$ is the vertical plane. In fact, if we replace, (\ref{vprimfun}) and (\ref{vsegfun}) in the mean curvature formula  (\ref{mean curvatue}), when $H=0$, we obtain that $\ddot{a}(t)=0$. That is, the minimal surface is a vertical plane.

Now we study the minimal graphs of $\mathcal{H}_{3},$ whose Gauss map have
rank $1.$ That is,%
\begin{equation}
f_{xx}f_{yy}-f_{xy}^{2}+\frac{1}{4}=0  \label{rank1}
\end{equation}

\begin{lemma}
\label{lemrul}Let $( x,y,f( x,y) ) $, with $(x,y)\in \Omega,$ be a minimal
graph in $\mathcal{H}_{3},$ which contains the origin, its normal at the
origin is $\eta (0) =\frac{1}{\sqrt{1+4k^{2}}}(0,-2k,1) $
and its Gauss map has rank 1. Furthermore assume that $%
f_{yy}(0) =0,$ then%
\[
f( x,y) =\left\{
\begin{array}{l}
\frac{xy}{2}+k\left[ \ln \left( y+\sqrt{1+y^{2}}\right) +y\sqrt{1+y^{2}}%
\right] \\
or \\
2ky-\frac{xy}{2}%
\end{array}%
\right.
\]
\end{lemma}

\begin{proof}
Since the unit normal at 0 is $\eta (0) =\frac{1}{\sqrt{1+4k^{2}}%
}\left( 0,-2k,1\right) $ we have, using $\left( \ref{normal}\right),$ that $%
f_{x}\left( 0\right) =0$ and $f_{y}\left( 0\right) =2k.$ On the other hand
 the Gauss map of such surface has rank 1, then, using $\left( \ref%
{rank1}\right),$ we have %
\[
f_{xx}\left( 0\right) f_{yy}\left( 0\right) -f_{xy}^{2}\left( 0\right) +%
\frac{1}{4}=0
\]%
 Since $f$ satisfies the minimal graph equation (\ref{minec}), we obtain
 $$(1+4k^{2})f_{xx}(0)+f_{yy}(0)=0.$$
 From the above two equations and the hypothesis, $f_{yy}(0)=0$,  we obtain
  $f_{xy}\left( 0\right) =\pm\frac{1}{2}$ and $f_{xx}\left( 0\right) =0.$

Recalling that $f$ is an analytic function, we can write its Taylor expansion in the form $f\left( x,y\right)
=2ky\pm \frac{xy}{2}+\Psi \left( x,y\right) .$ Substituting into $\left( \ref%
{rank1}\right) $ we obtain%
\[
\Psi _{xx} \Psi _{yy} -\Psi _{xy}^{2} =\pm \Psi _{xy}
\]%
Let $n$ be the minimal order of $\Psi \left( x,y\right).$ We claim that the
terms of minimal order of $\Psi $ do not appear mixed. In fact, assuming
that this is not the case, the minimal order of $\Psi _{xy}$ and $\Psi
_{xx}\Psi _{yy}-\Psi _{xy}^{2}$ are $\left( n-2\right) $ and $2\left(
n-2\right) $ respectively. This is a contradiction and proves our claim.

We shall now compute the third partial derivatives of $f$ at $0.$ To do this
we differentiate, with respect to $x$ and $y$, the minimal graph equation $%
\left( \ref{minec}\right) ,$ the equation $\left( \ref{rank1}\right) $ and
evaluate both at $0.$ We obtain the following two cases:

1.  If $f_{xy}\left( 0\right) =\frac{1}{2},$ we have that $f_{yyy}\left(
0\right) =2k$ and the others third derivatives of $f$ are zero. Then the
Taylor expansion of $f$ has the form:%
\[
f(x,y)=2ky+\frac{xy}{2}+\psi (y)+ax^{n}+by^{n}+\tilde{\Psi}(x,y),
\]%
where $\psi(y) $ is a polynomial such that $3\leq deg \psi(y)\leq (n-1)$
 with $n\geq 4.$ Whence%
\[
\begin{array}{rcl}
f_{x} & = & \frac{y}{2}+anx^{n-1}+\tilde{\Psi}_{x}, \\
f_{y} & = & 2k+\frac{x}{2}+\psi _{y}+bny^{n-1}+\tilde{\Psi}_{y}, \\
f_{xx} & = & an(n-1)x^{n-2}+\tilde{\Psi}_{xx}, \\
f_{xy} & = & \frac{1}{2}+\tilde{\Psi}_{xy}, \\
f_{yy} & = & \psi _{yy}+bn(n-1)y^{n-2}+\tilde{\Psi}_{yy}.%
\end{array}%
\]%
Substituting into $(\ref{minec}) $

\noindent
$
 (an(n-1)x^{n-2}+\tilde{\Psi}_{xx})[1+(2k+\psi _{y}+bny^{n-1}+\tilde{\Psi}%
_{y})^{2}]\\-2(\frac{1}{2}+\tilde{\Psi}_{xy})(2k+\psi _{y}+bny^{n-1}+\tilde{%
\Psi}_{y})(y+anx^{n-1}+\tilde{\Psi}_{x})\\+
(\psi _{yy}+bn(n-1)y^{n-2}+\tilde{\Psi}_{yy})[1+(y+anx^{n-1}+\tilde{\Psi}%
_{x})^{2}]=0.
$

\noindent If we analyze the coefficient of the term $x^{n-2}$ we obtain that%
\[
(1+4k^{2})an(n-1)x^{n-2}=0
\]%
Hence $a=0$ that is, $f\left( x,y\right) =\frac{xy}{2}+g\left( y\right) .$
We conclude that such a surface is invariant under translation of type $%
L_{\left( b,0,0\right) },$ see \cite{figuer},  and therefore
\[
f(x,y)=\frac{xy}{2}+k[\ln (y+\sqrt{1+y^{2}})+y\sqrt{1+y^{2}}].
\]%
for some $k\in \mathbb{R}$

2.  If $f_{xy}(0) =-\frac{1}{2},$ we have that the third
partial derivatives of $f$, evaluated at the origin, are equal to zero;
then, the Taylor expansion of $f$ has the form: $f( x,y) =$ $2ky-%
\frac{xy}{2}+ax^{n}+by^{n}+\tilde{\Psi}(x,y),$ where $n\geq 4.$ Whence%
\[
\begin{array}{rcl}
f_{x} & = & -\frac{y}{2}+anx^{n-1}+\tilde{\Psi}_{x}, \\
f_{y} & = & 2k-\frac{x}{2}+bny^{n-1}+\tilde{\Psi}_{y}, \\
f_{xx} & = & an(n-1)x^{n-2}+\tilde{\Psi}_{xx}, \\
f_{xy} & = & \tilde{\Psi}_{xy}-\frac{1}{2}, \\
f_{yy} & = & bn(n-1)y^{n-2}+\tilde{\Psi}_{yy}.%
\end{array}%
\]%
Substituting into $(\ref{minec}) $ we obtain

$(an(n-1)x^{n-2}+\tilde{\Psi}_{xx})[1+(2k-x+bny^{n-1}+\tilde{\Psi}%
_{y})^{2}]-2(\tilde{\Psi}_{xy}-\frac{1}{2})(2k-x+bny^{n-1}+\tilde{\Psi}%
_{y})(anx^{n-1}+\tilde{\Psi}_{x})+(bn(n-1)y^{n-2}+
\tilde{\Psi}_{yy})[1+(anx^{n-1}+\tilde{\Psi}_{x})^{2}]=0.$

\noindent If we analyze the coefficients of $x^{n-2}$ and $y^{n-2},$ we conclude that $%
a=b=0.$ Therefore,
\[
f(x,y)=2ky-\frac{xy}{2}.
\]%
This conclude the proof.
\end{proof}

We shall prove now that every minimal graph with Gauss map of rank 1 must
be a ruled surface. Let $S$ be such a surface, parameterized as a graph of a
differentiable function $f,$ with $f\left( 0,0\right) =0.$ Since $S$ has
rank 1, there exists a curve in $S,$ passing through the origin, such that
the unit normal field along this curve is constant. We indicate this curve
by $\Gamma ,$ where%
\[
\Gamma (t):\left\{
\begin{array}{ccl}
x(t) & = & t \\
y(t) & = & \alpha (t) \\
z(t) & = & f(t,\alpha (t))%
\end{array}%
\right. ,\,\,\,t\in (-\epsilon ,\epsilon ),
\]%
where $\alpha(0)=0$. We can assume that the normal field at $0$ is $\eta(0)=\frac{1}{\sqrt{1+4k^{2}}}(0,-2k,1).$
 Then, along the curve $\Gamma( t), $ and using $( \ref{normal}) $ we obtain,%
\begin{equation}
\begin{array}{rll}
f_{x}(t)+\frac{\alpha (t)}{2} & = & 0 \\
f_{y}(t)-\frac{t}{2} & = & 2k%
\end{array}%
\,,\,\,\,t\in (-\epsilon ,\epsilon ).  \label{ecucur}
\end{equation}%
Whence,
\begin{equation}
\begin{array}{rll}
f_{xx}(t)+\alpha ^{\prime }f_{xy}(t)+\frac{\alpha ^{\prime }(t)}{2} & = & 0
\\
f_{yx}(t)+\alpha ^{\prime }f_{yy}(t)-\frac{1}{2} & = & 0.%
\end{array}
\label{eccur2}
\end{equation}%
We need also the second and third partial derivatives of $f$ evaluated along
the curve $\Gamma \left( t\right) .$ From $\left( \ref{ecucur}\right) $ and
equations $\left( \ref{minec}\right) $ and $\left( \ref{rank1}\right) ,$ we
obtain the followings expressions for the partial derivatives of $f:$%

\[
\begin{array}{ccl}
f_{xx}(t) & = & \dfrac{\sin \theta }{2\sqrt{1+4k^{2}}} \\
f_{xy}(t) & = & \dfrac{\cos \theta }{2} \\
f_{yy}(t) & = & -\dfrac{\sqrt{1+4k^{2}}\sin \theta }{2}%
\end{array}%
\]%
Using this partial derivatives and the second equation of (\ref{eccur2}),  we obtain
\begin{equation} \label{alfa1}
\alpha'(t)=\frac{cos\theta-1}{sin\theta\sqrt{1+4k^{2}}}.
\end{equation}
Now, by differentiating the equations $\left( \ref{minec}\right) $ and $\left( %
\ref{rank1}\right) ,$ with respect to $x$ and $y,$ and evaluating at $\left(
t,\alpha \left( t\right) \right) ,$ we obtain the following system%
\[
\left(
\begin{array}{cccc}
\frac{-\sqrt{1+4k^{2}}\sin \theta }{2} & -\cos \theta  & \frac{\sin \theta }{%
2\sqrt{1+4k^{2}}} & 0 \\
0 & \frac{-\sqrt{1+4k^{2}}\sin \theta }{2} & -\cos \theta  & \frac{\sin
\theta }{2\sqrt{1+4k^{2}}} \\
(1+4k^{2}) & 0 & 1 & 0 \\
0 & (1+4k^{2}) & 0 & 1%
\end{array}%
\right) \left(
\begin{array}{c}
f_{xxx}(t) \\
f_{xxy}(t) \\
f_{xyy}(t) \\
f_{yyy}(t)%
\end{array}%
\right) =\left(
\begin{array}{c}
0 \\
0 \\
\frac{k\sin \theta }{\sqrt{1+4k^{2}}} \\
k(1+\cos \theta )%
\end{array}%
\right)
\]%
Solving this system, we obtain the third partial derivatives of
$f$ along the
curve $\Gamma \left( t\right) :$%
\begin{equation}
\begin{array}{ccl}
f_{xxx}(t) & = & \dfrac{k\sin \theta (1-\cos \theta )}{2(1+4k^{2})^{3/2}} \\
&  &  \\
f_{xxy}(t) & = & \dfrac{k\sin {}^{2}\theta }{2(1+4k^{2})} \\
&  &  \\
f_{yyy}(t) & = & \dfrac{k(1+\cos \theta )^{2}}{2} \\
&  &  \\
f_{yyx}(t) & = & \dfrac{k\sin \theta (1+\cos \theta )}{2\sqrt{1+4k^{2}}}.%
\end{array}
\label{terceras}
\end{equation}%
We are now ready to prove

\begin{theorem}
\label{classi}If $\left( x,y,f\left( x,y\right) \right) $ with $(x,y)\in\Omega\subset \mathbb{R}^{2}$ is a minimal
graph such that its normal at the origin is $\eta \left( 0\right) =\frac{1}{%
\sqrt{1+4k^{2}}}(0,-2k,1)$ and its Gauss map has rank 1, then it is a ruled
surface.
\end{theorem}

\begin{proof}
We shall show that the curve $\Gamma ,$ defined above, must be a straight
line. To do this , we differentiate the third component of $\Gamma ,$ that
is,%
\[
\frac{dz}{dt}=f_{x}(t)+\alpha ^{\prime }f_{y}(t)
\]%
From this expression and $\left( \ref{eccur2}\right) ,$ we have that $\dfrac{%
d^{2}z}{dt^{2}}=$\- $\alpha ^{\prime \prime }f_{y}(t).$ By differentiating
the second equation of $\left( \ref{eccur2}\right) $ with respect to $t,$ we
obtain%
\begin{equation}
f_{yxx}(t)+2\alpha ^{\prime }f_{yyx}(t)+(\alpha ^{\prime
})^{2}f_{yyy}(t)+\alpha ^{\prime \prime }f_{yy}(t)=0.  \label{ecterc}
\end{equation}%
Then we have two cases:

\begin{enumerate}
\item If $f_{yy}\left( 0\right) =0,$ using the proposition $\left( \ref%
{lemrul}\right) ,$ we conclude that such a surface is a ruled surface

\item If $f_{yy}\left( 0\right) \neq 0$ then, $f_{yy}\left( t\right) \neq 0,$
$t\in \left( -\varepsilon ,\varepsilon \right) $ By replacing ( \ref%
{terceras}) and ( \ref{alfa1}) in ( \ref{ecterc}) , we obtain that $\alpha ^{\prime \prime
}f_{yy}\left( t\right) =0,$ $t\in ( -\varepsilon ,\varepsilon) .$
Then $\alpha \left( t\right) =bt$ and therefore, from ( \ref{ecucur}),
 we conclude that $f( t,bt) =2kbt$ . This completes the
proof.
\end{enumerate}
\end{proof}

The following classification result for ruled minimal surfaces in $\mathcal{H%
}_{3}$ was proved by M. Bekkar and T. Sari, see \cite{bekkar}.

\begin{theorem}
The ruled minimal surfaces of $\mathcal{H}_{3},$ up isometry, are:

\begin{enumerate}
\item The plane;

\item The hyperbolic paraboloid;

\item The helicoid parameterized by%
\[
\left\{
\begin{array}{ccl}
x(t,s) & = & s\sin t \\
y(t,s) & = & s\cos t \\
z(t,s) & = & \rho t%
\end{array}%
,\,\,\,\,\,\,\,\rho \in I\!\!R-\{0\}.\right.
\]

\item The surface of equation%
\[
z=\frac{xy}{2}-\frac{\lambda }{2}\left[ y\sqrt{1+y^{2}}+\log \left( y+\sqrt{%
1+y^{2}}\right) \right] ,\,\,\,\,\,\lambda \in I\!\!R-\{0\}.
\]

\item The surfaces which are locally the graph of the function $z=\frac{y}{2}%
\left( R\left( x\right) +x\right) ,$ where $R$ is a solution of the
differential equation%
\[
R^{\prime \prime }\,\,\left( 4+R^{2}\right) -2R\left( R^{\prime }+1\,\right)
\left( R^{\prime }+2\,\right) =0.
\]

\item The surfaces which are locally parameterized by
\[
\left\{
\begin{array}{ccl}
x(t,s) & = & t+su(t) \\
y(t,s) & = & s \\
z(t,s) & = & a(t)-\frac{st}{2}%
\end{array}%
\right.
\]%
where $u$ and $a$ are solutions of the system%
\begin{equation}
\left\{
\begin{array}{lcc}
\left( 1+u^{2}+t^{2}\right) u^{\prime \prime }-\left( 1+2u^{\prime
}a^{\prime }\right) tu^{\prime } & = & 0 \\
\left( 1+u^{2}+t^{2}\right) a^{\prime \prime }-\left( 1+2u^{\prime
}a^{\prime }\right) (ta^{\prime }-u) & = & 0.%
\end{array}%
\right.  \label{item6}
\end{equation}
\end{enumerate}
\end{theorem}

The above theorem, together with theorem $\left( \ref{classi}\right) ,$
give us the classification for  minimal surfaces in $\mathcal{H}_{3}$ with
Gauss map of rank 1. In fact, to do this we determine, among the surfaces given by
the above classification, which ones have the rank equal to 1.

It is not difficult to compute,using $\left( \ref{rank1}\right) $ that the
rank of the Gauss map of the surfaces of items 1 and 3 of the above theorem
have rank different from 1 and the surfaces of items 2 and 4 have rank 1.

We shall now study the surface of the item 5, In this case , we have that $%
f\left( x,y\right) =\frac{y}{2}\left( R\left( x\right) +x\right) .$ Hence,%
\[
f_{xx}=R^{\prime \prime };\,\,\,\,\,f_{xy}=\frac{1}{2}\left( R\,^{\prime
}+1\,\right) ;\,\,\,\,\,f_{yy}=0.
\]%
Substituting into $\left( \ref{rank1}\right) $ , we obtain $\left( R^{\prime
}+1\right) ^{2}=1.$ Solving this differential equation we obtain the
following solution%
\[
f\left( x,y\right) =\frac{y}{2}\left( a-x\right).
\]%
The isometry $L_{\left( a,0,0\right) }$ takes this surface to the parabolic
hyperboloid $z=\frac{xy}{2}.$

In the case of the surface of item 6, we parameterized this surface as a
graph of a differential function $f,$ where%
\[
f\left( x,y\right) =a\left( t\left( x,y\right) \right) -\frac{1}{2}yt\left(
x,y\right)
\]%
Whence
\[
\left\{
\begin{array}{l}
t_{x}=\frac{1}{1+yu^{\prime }};\,\,\,\,\,t_{y}=\frac{-u}{1+yu^{\prime }}%
\,\,\, \\
f_{xx}=\,a^{\prime \prime }t_{x}^{2}+\left( a^{\prime }+\frac{y}{2}\,\right)
t_{xx} \\
f_{xy}=\frac{1}{2}t_{x}+a^{\prime \prime }t_{x}t_{y}+\left( a^{\prime }\,+%
\frac{y}{2}\right) t_{xy} \\
f_{yy}=a^{\prime \prime }t_{y}^{2}+t_{y}+\left( a^{\prime }+\frac{y}{2}%
\right) t_{yy}%
\end{array}%
\right.
\]%
Substituting into $\left( \ref{rank1}\right) ,$ we obtain that $t_{x}^{2}%
\left[ \frac{1}{2}-\left( a^{\prime }+\frac{y}{2}\right) u^{\prime }t_{x}%
\right] ^{2}=\frac{1}{4}.$ Whence%
\[
\left( u^{\prime }\right) ^{2}y^{2}+2u^{\prime }y+2u^{\prime }a^{\prime }=0.
\]%
Notice that this is a polynomial of second degree with respect to $y$ and
its coefficients depend only on $t$, so, there must be zero, that is, $u^{\prime }=0.$ Hence $u$ is
constant. By mean a rotation about the $z-axis$ we may take $u=0.$ By
replacing in the second equation of  $\left( \ref{item6}\right) ,$ we obtain%
\[
\left( 1+t^{2}\right) a^{\prime \prime }-ta^{\prime }=0.
\]%
The general solution of this equation is
\[
a\left( t\right) =\frac{\lambda }{2}\left[ t\sqrt{1+t^{2}}+\ln \left( t+%
\sqrt{1+t^{2}}\right) \right] +\mu .
\]%
where $\lambda ,\mu \in \mathbb{R}.$ This minimal surface may be expressed
as the graph of the function%
\[
f\left( x,y\right) =\frac{xy}{2}-\frac{\lambda }{2}\left[ y\sqrt{1+y^{2}}%
+\ln \left( y+\sqrt{1+y^{2}}\right) \right]
\]

Therefore we have the following classification for minimal surfaces in $%
\mathcal{H}_{3}.$

\begin{theorem}
The minimal graphas in $\mathcal{H}_{3}$ with Gauss map of rank 1, are
\[
f\left( x,y\right) =\frac{xy}{2}-\frac{k}{2}\left[ y\sqrt{1+y^{2}}+\ln
\left( y+\sqrt{1+y^{2}}\right) \right]
\]%
where $k\in \mathbb{R}$
\end{theorem}

To conclude this section we shall present some results about complete
minimal graphs in $\mathcal{H}_{3}$. Firstly, we present one directly consequence
of the minimal graph equation.

\begin{proposition}
\label{acotado}If $f\left( x,y\right) $ is a function that satisfies $\left( %
\ref{minec}\right) ,$ then%
\[
f_{xx}f_{yy}-f_{xy}^{2}\leq 0,\ \ \ \ \forall \left( x,y\right).
\]
\end{proposition}

\begin{proof}
Let $a=1+\left( f_{x}+\frac{y}{2}\right) ^{2},$ $b=\left( f_{x}+\frac{y}{2}%
\right) \left( f_{y}-\frac{x}{2}\right) $ and $c=1+\left( f_{y}-\frac{x}{2}%
\right) ^{2}.$ Then using equation $\left( \ref{minec}\right) ,$ we obtain%
\[
f_{xx}f_{yy}-f_{xy}^{2}=-\frac{1}{a}\left(
af_{xy}^{2}+2bf_{yy}f_{xy}+cf_{yy}^{2}\right)
\]%
Since $ac-b^{2}=1+\left( f_{x}+\frac{y}{2}\right) ^{2}+\left( f_{y}-\frac{x}{%
2}\right) ^{2}>0,$ the result follows
\end{proof}

Now we recall the following theorem of Bernstein, see \cite{hoph}and
\cite{mick}.

\begin{theorem}
Let $f\left( x,y\right) $ be a real-valued function which satisfies the
following conditions

\begin{enumerate}
\item $f\left( x,y\right) \in C^{2}\left( \mathbb{R}^{2}\right) .$

\item $f_{xx}f_{yy}-f_{xy}^{2}\leq 0,$ $f_{xx}f_{yy}-f_{xy}^{2}\neq 0.$
\end{enumerate}

Then $f\left( x,y\right) $ is not bounded.
\end{theorem}

It follows from the above theorem that a complete minimal graph in $\mathcal{%
H}_{3}$ cannot be bounded. More precisely.

\begin{proposition}
Let $\left( x,y,f\left( x,y\right) \right) $ be a minimal graph in $\mathcal{%
H}_{3},$ which is defined in the entire $xy-plane$ and $\left\vert f\left(
x,y\right) \right\vert \leq k,\ \ \forall \left( x,y\right) .$ Then $f$ is
constant.
\end{proposition}

\begin{proof}
Since $f$ is bounded we have, using the Bernstein's theorem that%
\[
f_{xx}f_{yy}-f_{xy}^{2}\equiv 0
\]%
It follows from the proof of the proposition $\left( \ref{acotado}\right) $
\[
af_{xy}^{2}+2bf_{yy}f_{xy}+cf_{yy}^{2}=0
\]%
where $a,b$ and $c$ are as in such proof. Then $f_{yy}=f_{xy}=0$ and,
substituting in minimal graph equation, we obtain that $f_{xx}=0.$
Therefore, $f$ is constant.
\end{proof}

A consequence of the above proof is that when a rank of the Gauss map of a complete minimal graph
is equal to $1/4w^{2}$ where $w$ is like (\ref{normal}), the surface must be a plane.

Finally, we must mentioned that I. Fernadez and P. Mira, gave a classification of the entire minimal vertical graphs in
$\mathcal{H}_{3}$ in terms of the Abresh-Rosenberg holomorphic differential for minimal surfaces in
$\mathcal{H}_{3}$, see \cite{femi}, but it is interesting to study the image of minimal surface when its Gauss map
 has rank 2, see, for example, \cite{bdani}.

\end{document}